\documentclass[]{amsart}

\usepackage{color}
\usepackage[]{hyperref}

\usepackage{natbib}

\usepackage{amssymb}

\newcommand{\Real}{\mathbb R}

\newcommand{\sfE}{\mbox{\sf E}}
\newcommand{\I}{\mbox{\sc I}}
\newcommand{\F}{{\mathcal F}}
\newcommand{\T}{{\mathcal T}}

\newcommand{\bP}{\mathbf{P}}

\newtheorem{theorem}{Theorem}[section]
\newtheorem{lemma}[theorem]{Lemma}
\newtheorem{corollary}[theorem]{Corollary}

\begin{document}
\title[On optimal stopping of a risk process]{Optimal stopping of a risk process when claims are covered immediately}
\thanks{This paper is a preliminary version, and the final form will be published elsewhere.}
\author[B.K. Muciek]{Bogdan K. Muciek}
\address{Institute of Mathematics and Computer Science, %
Wroc\l{}aw University of Technology, Wybrze\.{z}e Wyspia\'{n}skiego 27, %
Wroc\l{}aw, Poland}
\curraddr{AIG Credit SA, ul. Strzegomska 42c, 53-611 Wroclaw, Poland}
\email{Bogdan.Muciek@aigcredit.pl}
\urladdr{http://www.im.pwr.wroc.pl/~muciek}
\thanks{The research was supported by KBN grant no  2 P03A 021 22 (350228).}
\author[K.J. Szajowski]{Krzysztof J. Szajowski}%
\address[Corresponding author]{Institute of Mathematics and Computer Science, Wroc\l{}aw University of Technology, Wybrze\.{z}e Wyspia\'{n}skiego 27, Wroc\l{}aw, Poland}
\email{Krzysztof.Szajowski@pwr.wroc.pl}
\urladdr{http://neyman.im.pwr.wroc.pl/\~{}szajow}

\begin{abstract}
The optimal stopping problem for the risk process with interests rates and when claims are covered
immediately is considered. An insurance company receives premiums and pays out claims which have occured according to a renewal process and which have been recognized by them. The capital of the company is invested at interest rate $\alpha\in\Re^{+}$, the size of claims increase at rate $\beta\in\Re^{+}$ according to inflation process. The immediate payment of claims decreases the company investment by rate $\alpha_1$. The aim is to find the stopping time which maximizes the capital of the company. The improvement to the known models by taking into account different scheme of claims payment and the possibility of rejection of the request by the insurance company is made. It leads to essentially new risk process and the solution of optimal stopping problem is different.
\end{abstract}



\thanks{Subject Category and Insurance Branch Category: IE43 IM51 IM10}

\maketitle

\section{Introduction}
The following problem in collective risk theory (see \cite{rolschschteu98}) is considered.
An insurance company, endowed with an initial capital $a>0$, receives premiums and pays out claims that occur
according to a renewal process $\{N(t),t\ge 0\}$, where $N(t)$ is the number of losses up till time $t$.
The initial capital of the insurance company and received premiums are invested at a constant rate of return
$\alpha\in\Re^{+}$. Let $T_0=0$ and let $T_i$ denotes the time of the $i$-th loss, then random variables
$\zeta_i=T_i-T_{i-1}$ are independent and identically distributed (i.i.d.)\ with cumulative  distribution
function (cdf) $F$, such that $F(0)=0$. Let $X_1,X_2,\ldots$ be a sequence of i.i.d.\ random variables
independent of $\{\zeta_i\}$, with cdf $H$ with $H(0)=0$. The sequence $\{X_i\}_{i=1}^\infty$ represents values
of successive claims. Usually the costs of damages elimination increase. It is modelled by the rate $\beta\in\Re^{+}$.
If a claim appears at moment $T_n$, then the company have to pay  $X_n e^{\beta T_n}$. This amount of money
decreases the company investment by rate $\alpha_1$ and, as a consequence of that, at the end of investment period
$t$ the claim at $T_n$ decreases the capital by $X_ne^{\beta T_n}e^{\alpha_1(t-T_n)}$. Although it may seem somewhat
surprising at first glance, claims of size zero arise in some insurance contexts (see \cite{panwil92:insurance}).
If a company records all claims as they are presented to the company and some claims are resisted, refused or
a complete recovery of losses is made from another insurer, the net cost of the claim is zero.
This effect is modelled by additional sequence of i.i.d random variables $\{\epsilon_i\}_{i=1}^\infty$,
independent of claims size process and the process of moments of claims. It is assumed that
$\bP\{\epsilon_n=1\}=p$ and $\bP\{\epsilon_n=0\}=1-p$.
The investigated process of capital assets of the insurance company is
\begin{equation}
U_t=a e^{\alpha t}+\int_0^tce^{\alpha(t-s)}ds-\sum_{n=0}^{N(t)}\epsilon_n X_ne^{\beta T_n}e^{\alpha_1(t-T_n)},
\label{fer1}
\end{equation}
where $a>0$ is the initial capital, $c>0$ is a constant rate of income from the insurance
premiums, $X_0=0$ and $N(0)=0$. The form of capital assets (\ref{fer1}) reduces to
\begin{eqnarray} \label{fer1a}
U_t&=&ae^{\alpha t}+ce^{\alpha t}\frac{1-e^{-\alpha t}}{\alpha}-e^{\alpha_1 t}\sum_{n=0}^{N(t)}%
\epsilon_n X_ne^{\beta_1 T_n}
\end{eqnarray}
where $\beta_1=\beta-\alpha_1$. Let $g(u,t)=g_1(u)\I_{\{t\geq 0\}}$, where $g_1$ is a utility function.
The return at time $t$ is $\{Z(t), t\ge 0\}$ and it is given by
\begin{equation}\label{fer2}
Z(t)=g(U_t,t_0-t)\prod_{j=0}^{N(t)}\I_{\{U_{T_j}>0\}}=g(U_t)\I_{\{U_s>0,s\leq t\}}
\end{equation}
The optimal stopping problem for the process $Z(t)$ is investigated. The model with $\alpha=\beta=0$
have been considered by \cite{fersie97:risk}.
\cite{jen97:risk} 
investigated a similar model with a claim process modulated by periodic Markovian processes but without
care for time value of money, formulated in (\ref{fer1a}). When the claims are paid from
the capital of the company, it can be assumed that $\alpha=\alpha_1$.
\cite{muc02:interest} 
investigated the model given by (\ref{fer1}) with $\alpha_1=0$ which described the case when the claims
were paid at the end of the investing period. The improvement introduced
here, which takes into account the consequence of the immediate payment of
claims, change the considered risk process essentially. The model admitted
will have an impact on the form of the strong generator for the process
(\ref{fer1a}) 
as well as on the form of the dynamic programming equations, which are the tools for describing the
solution of the optimal stopping problem for (\ref{fer2}).

The considered process $Z(t)$ is the piecewise-deterministic process. The methods of solving the optimal stopping
problem for such processes can by found in papers by
\cite{bosgou93:semi} and \cite{jen97:risk,jenhsu93:point}
and the monography by \cite{davis93:MModels}. \cite{muc02:interest} has solved
the optimal stopping problem for process (\ref{fer2}) with $\alpha_1=0$ which is not direct consequence
of the optimal stopping problem solution for model (\ref{fer1}) with $\alpha_1\neq 0$.

The organization of the paper is following. In the next section the optimal stopping problem for the process (\ref{fer2})
is formulated. The case of the optimal stopping up to the fixed number of claims is the subject of investigation in
the section \ref{fixedclaim}. The solution of the problem for the infinite number of claims is given in the section
\ref{infcase}.

\section{The optimization problem}

In this section we define an optimization problem for the model introduced in the
previous section. This optimization problem will be solved in the next section.

Let
$
{\mathcal F}(t)=\sigma(U_s,s\le t)=\sigma(X_1,\epsilon_1,T_1,\ldots,X_{N(t)},\epsilon_{N(t)},T_{N(t)})
$
be the $\sigma$-field generated by all the events up to time $t\ge 0$.
Let $\mathcal T$ be the
set of all stopping times with respect to the family $\{{\mathcal F}(t),t\ge 0\}$. Furthermore, for fixed $K$ and
for
$n=0,1,\ldots,k<K$ let ${\mathcal T}_{n,K}$ denote the subset of $\mathcal T$, such that
$\tau\in{\mathcal T}_{n,K}$ if and only if $T_n\le\tau\le T_K \quad\mbox{ a.s.}
$

Let $\F_n=\F(T_n)$. The essence of the
considerations in the next chapter will be to find the optimal stopping time $\tau^*_K$, such that
$\sfE Z(\tau^*_K)=\sup\{\sfE Z(\tau):\:\tau\in\T_{0,K}\}$.
In order to find the optimal stopping time $\tau^*_K$, we first consider optimal stopping
times $\tau^*_{n,K}$, such that
\begin{equation}
\sfE(Z(\tau^*_{n,K})|\F_n)=\mbox{ess}\sup\{\sfE(Z(\tau)|\F_n):\:\tau\in\T_{n,K}\}
\label{taunk}
\end{equation}
and using backward induction as in dynamic programming, we will obtain
$\tau^*_{K}=\tau^*_{0,K}$

After finding the optimal stopping time $\tau^*_K$ for fixed $K$ we will deal with unlimited number of claims
and the aim will be to find the optimal stopping time $\tau^{*}$, such that
\begin{equation}
\sfE Z(\tau^*)=\sup\{\sfE Z(\tau):\:\tau\in\T\}\label{fer6a}
\end{equation}
is fulfilled. It will be shown that $\tau^*$ can be defined as the limit of
the finite horizon optimal stopping times. Such a stopping time in an insurance company management can be used
as the best moment to recalculate premium rate.

\section{\label{fixedclaim}Case with fixed number of claims}

In this section we find the form of optimal stopping time in the finite horizon case, which means the optimal
stopping time in the class $\T_{0,K}$, where $K$ is finite and fixed (the number of claims is fixed, but the
time of the $K$th claim, ie.\ time horizon, remains non-deterministic). This is a technical assumption which
allows calculations for finite number of claims and generalize this result to the infinite number of claims.
First we present dynamic programming equations satisfying
$
\Gamma_{n,K}=\mbox{ess}\sup\{\sfE(Z(\tau)|\F_n):\tau\in\T_{n,K}\},\qquad n=K,K-1,\ldots,1.
$
Then in Corollary \ref{fertw2} we find optimal stopping times
$\tau^*_{n,K}$ and $\tau^*_K$ and optimal mean values of return related to them.

The following representation lemma (see for example \cite{davis76:representation}) 
plays the crucial role in consequent considerations:
\begin{lemma}\label{ferlem1}
If $\tau\in\T_{n,K}$, there exists a positive, $\F_n$-measurable random variable $\xi$, such
that $\tau\wedge T_{n+1}=(T_n+\xi)\wedge T_{n+1} \mbox{ a.s.}$
\end{lemma}
Let $\mu_0=1$ and $\mu_n=\prod_{j=1}^n\I_{\{U_{T_j}>0\}}$. Then
$
\Gamma_{K,K}=Z(T_K)=g(U_{T_K},t_0-T_K)\mu_K.
$
Note that the sum of claims from (\ref{fer1a}) can be expressed as
\begin{equation}\label{suma}
\sum_{n=0}^{N(t)}\epsilon_n X_ne^{\beta_1 T_n}=
\left(a e^{\alpha t}+\frac{c}{\alpha}(e^{\alpha t}-1)-U_t\right)e^{-\alpha_1 t}
\end{equation}
Let us define for $\xi>0$ such that there is no jump between $t$ and $t+\xi$
\begin{eqnarray}\label{d}
d(t,\xi,U_t)&=&U_{t+\xi}-U_t=e^{\alpha t}\left(a+\frac{c}{\alpha}\right)\left(e^{\alpha \xi}-e^{\alpha_1 \xi}\right)\\
\nonumber &&+\frac{c}{\alpha}\left(e^{\alpha_1\xi}-1\right)+\left(e^{\alpha_1\xi}-1\right)U_t,
\end{eqnarray}
then we have
\begin{equation}
\mu_K=\mu_{K-1}\I_{\{U_{T_{K-1}}+d(T_{K-1},\zeta_K,U_{T_{K-1}})-\epsilon_K X_Ke^{\beta(T_{K-1}+\zeta_K)}>0\}}.
\label{miu}
\end{equation}

Similarly as in \cite{muc02:interest}, Theorem 1, from (\ref{suma}) and from (\ref{d}) we get the following
dynamic programming equations:
\begin{description}
\item[(i)] For $n=K-1,K-2,\ldots,0$,
\begin{eqnarray*}
\Gamma_{n,K}&=&\mbox{\rm ess}\sup\big\{\mu_n\bar{F}(\xi)g\left(U_{T_n}+d(T_n,\xi,U_{T_n}),t_0-T_n-\xi\right)\\
&&+\sfE(\I_{\{\xi\ge \zeta_{n+1}\}}\Gamma_{n+1,K}|\F_n):\:\xi\ge 0
\mbox{ is }\F_n\mbox{-measurable}\big\}\quad \mbox{a.s.},
\end{eqnarray*}
where $\bar{F}=1-F$ is the survival function.

\item[(ii)] For $n=K,K-1,\ldots,0$,
$\Gamma_{n,K}=\mu_n\gamma_{K-n}(U_{T_n},T_n)\quad\mbox{a.s.}$,
where the sequence of functions $\{\gamma_j(u,t),u\in\Real,t\ge 0\}$, using (\ref{d}), (\ref{suma}) and (\ref{miu})
is defined as follows
{\jot0cm
\begin{eqnarray*}
&&\gamma_0(u,t)=g(u,t_0-t),\\
&&\gamma_j(u,t)=\sup_{r\ge 0}\bigg[\bar{F}(r)g(u+d(t,r,u),t_0-t-r)\\
\nonumber&&\mbox{$\;$}
+p\int_0^rdF(s)\int_0^{e^{-\beta(t+s)}(u+d(t,s,u))}\gamma_{j-1}\left(u+d(t,s,u)-xe^{\beta(t+s)},t+s\right)dH(x)\\
\mbox{$\;$}&&+(1-p)\int_0^r \gamma_{j-1}\left(u+d(t,s,u),t+s\right)H(e^{-\beta(t+s)}(u+d(t,s,u)))dF(s)\bigg]
\end{eqnarray*}}
$j=1,2,\ldots$
\end{description}

The above equations differ from the ones in Theorem 1 in \cite{muc02:interest} as a result of a different form of the 
capital assets process $U_t$.

The next step is to find the optimal stopping time $\tau^*_K$. To do this we should analyze
the properties of the sequence of functions $\{\gamma_n, n\ge0\}$.
Let $B=B[(-\infty,+\infty)\times[0,+\infty)]$ be the space of all bounded and continuous
functions with the norm $||\delta||=\sup_{u,t}|\delta(u,t)|$ and let
$
B^0=\{\delta:\:\delta(u,t)=\delta_1(u,t)\I_{\{t\le t_0\}}\mbox{ and }\delta_1\in B\}.
$
One should notice that the functions $\{\gamma_n, n\ge0\}$ are
included in $B^0$.
For each $\delta\in B^0$ and any $u\in\Real$, $t,r\ge0$ let
\begin{eqnarray*}
\phi_\delta(r,u,t)&=&\bar{F}(r)g(u+d(t,r,u),t_0-t-r)+\label{phideltaeq}\\
\nonumber&&+(1-p)\int_0^r \delta(u+d(t,s,u),t+s)H(e^{-\beta(t+s)}(u+d(t,s,u)))dF(s)\\
\nonumber&&+p\int_0^rdF(s)\int_0^{e^{-\beta(t+s)}(u+d(t,s,u))}\delta(u+d(t,s,u)-xe^{\beta(t+s)},t+s)dH(x).
\end{eqnarray*}

From the properties of the cumulative distribution function $F$ we know that
$\phi_\delta(r,u,t)$ has at most a countable number of points of discontinuity according to
$r$ and is continuous according to $(u,t)$ in the case of $g_1(\cdot)$ being continuous and
$t\not=t_0-r$. Therefore, for further considerations we assume that the function $g_1(\cdot)$ is bounded and continuous.

For each $\delta\in B^0$ let
\begin{equation}
(\Phi\delta)(u,t)=\sup_{r\ge 0}\{\phi_\delta(r,u,t)\}.
\label{fer18}
\end{equation}

\begin{lemma}
For each $\delta\in B^0$ we have
$$(\Phi\delta)(u,t)=\max_{0\le r\le t_0-t}\{\phi_\delta(r,u,t)\}\in B^0$$
and there exists a function $r_\delta(u,t)$ such that
$(\Phi\delta)(u,t)=\phi_\delta(r_\delta(u,t),u,t)$. \label{ferlem2}
\end{lemma}

In subsequent considerations more properties of $\Phi$ will be presented.

For $i=1,2,\ldots$ and $u\in\Real$, $t\ge 0$, $\gamma_i(u,t)$ may be expressed as follows
$$
\gamma_i(u,t)=\left\{
\begin{array}{ll}
(\Phi\gamma_{i-1})(u,t)&\mbox{ if }u\ge 0 \mbox{ and }t\le t_0,\\
0&\mbox{ otherwise,}
\end{array}\right.
$$
and from Lemma \ref{ferlem2} there exist functions $r_{\gamma_{i-1}}$ such that
$$
\gamma_i(u,t)=\left\{
\begin{array}{ll}
\phi_{\gamma_{i-1}}(r_{\gamma_{i-1}},u,t)&\mbox{ if }u\ge 0 \mbox{ and }t\le t_0,\\
0&\mbox{ otherwise.}
\end{array}\right.
$$

To specify the form of the optimal stopping times $\tau^*_{n,K}$, we need to define the
following random variables $R^*_i=r_{\gamma_{K-i+1}}(U_{T_i},T_i)$ and
$\sigma_{n,K}=K\wedge\inf\{i\ge n:\:R^*_i<\zeta_{i+1}\}$.

Finally in Corollary \ref{fertw2} we present the form of the optimal stopping time.

\begin{corollary}\label{fertw2}
Let
$$
\tau^*_{n,K}=T_{\sigma_{n,K}}+R^*_{\sigma_{n,K}}\quad \mbox{ and }\quad\tau^*_K=\tau^*_{0,K},
$$
then for all $0\le n\le K$ the following hold
$$
\Gamma_{n,K}=\sfE(Z(\tau^*_{n,K})|\F_n)\mbox{ a.s.}\quad\mbox{ and }\quad
\Gamma_{0,K}=\sfE(Z(\tau^*_K))=\gamma_K(a,0),
$$
which means $\tau^*_{n,K}$ and $\tau^*_K$ are optimal stopping times in the classes $\T_{n,K}$
and $\T_{0,K}$ respectively.
\end{corollary}

\section{\label{infcase}Case with an infinite number of claims}
While $\mathcal T$ is the set of all stopping times with respect to the family
$\{{\mathcal F}(t),t\ge0\}$, we would like to maximize the mean return
(\ref{fer2}), i.e. to find the optimal stopping time $\tau^*$, such that

\begin{equation}
\sfE Z(\tau^*)=\sup\{\sfE Z(\tau):\:\tau\in\T\}\label{fer6}
\end{equation}
is fulfilled. It will be shown that $\tau^*$ can be defined as the limit of
the finite horizon optimal stopping times.

Let $\tau^*_n$ be an optimal stopping time taken from the set of stopping times which occur not earlier than
at $T_n$, the time of $n$-th claim, ie.
$\sfE Z(\tau^*_n|\F_n)=\sup\{\sfE Z(\tau|\F_n):\:\tau\in\T\cap\{\tau:\tau\ge T_n\}\}$.

The solution of this case will be based on the iteration of the operator $\Phi$ defined by (\ref{fer18}). By
an assumption that the interarrival time is greater than $t_0$ with non-zero probability it can proved that
the operator $\Phi$ is a contraction and it has a fixed point.

The essence of this section is contained in the following theorem. The proof is
based on the proof of the existence of optimal stopping times for semi-Markov processes
presented by \cite{bosgou93:semi}. 

\begin{theorem}
Assuming the utility function $g_1$ is differentiable and nondecreasing and $F$ has a density
function $f$ we have
\begin{itemize}
\item[\em(i)]
for $n=0,1,\ldots$, the limit $\tilde{\tau}_n:=\lim_{K\to\infty}\tau_{n,K}^*$ exists and
$\tilde{\tau}_n$ is an optimal stopping time in $\T\cap\{\tau:\;\tau\ge T_n\}$ ($\tau^{*}_{n,K}$ is a
solution of the case with finite number of claims defined by (\ref{taunk})),
\item[\em(ii)]
$\sfE[Z(\tilde{\tau}_n)|\F_n]=\mu_n\gamma(U_{T_n},T_n)$\quad\mbox{a.s.}
\end{itemize}
\end{theorem}

The optimality of $\tilde{\tau}_n$ may be proved in a similar way as in \cite{bosgou93:semi}. 
The details may be found in \cite{muc02:interest}. 
Let us denote $\eta_t=(t,U_t,Y_t,V_t)$, where  $Y_t=t-T_{N(t)}$, $V_t=\mu_{N(t)}$, $t\ge0$. For $\tilde{g}$
such that $Z(t)=\tilde{g}(\eta_t)$ the strong generator $A$ of $\eta_t$ on $\tilde{g}$ has the form
(see \cite{giksko71:processes2E} 
)
{\small
\begin{eqnarray*}
(A\tilde{g})(t,u,y,v)&=&\left\{\bigg[(\alpha a+c)e^{\alpha t}-(\alpha_1 a+c\frac{\alpha_1}{\alpha})e^{\alpha t}
+ (\frac{\alpha_1}{\alpha}c+\alpha_1 u)\bigg]g'_1(u)\right.\\
&&-\frac{f(y)}{\bar{F}(y)}\bigg[g_1(u)\bar{H}(ue^{-(\alpha_1 y+\beta (t-y))})\\
&&-p\int_0^{ue^{-(\alpha_1 y+\beta (t-y))}}g_1(u-xe^{\alpha_1 y+\beta (t-y)})dH(x)\\
&&+pg_1(u)\left.H(u e^{-(\alpha_1 y+\beta (t-y))})\bigg]\right\}v,
\end{eqnarray*}
}
where $t<t_0$, $y\ge0$ and $v\in\{0,1\}$.
Thus we get
{\small
\begin{eqnarray*}
(A\tilde{g})(\eta_s)
&=&\left\{\bigg[(\alpha a+c)e^{\alpha s}-(\alpha_1 a+c\frac{\alpha_1}{\alpha})e^{\alpha s}
+ (\frac{\alpha_1}{\alpha}c+\alpha_1 U_s)\bigg]g'_1(U_s)\right.\\
&+&\frac{f(s-T_{N(s)})}{\bar{F}(s-T_{N(s)})}
\bigg[
p\int_0^{U_s e^{-(\alpha_1 (s-T_{N(s)})+\beta T_{N(s)})}}
g_1(U_s-xe^{\alpha_1 (s-T_{N(s)})+\beta T_{N(s)}})dH(x)\\
&&-pg_1(U_s)H(U_s e^{-(\alpha_1 (s-T_{N(s)})+\beta T_{N(s)})})\\
&&-\left.\bar{H}(U_s e^{-(\alpha_1 (s-T_{N(s)})+\beta T_{N(s)})})g_1(U_s)\bigg]\right\}
\mu_{N(s)}.
\end{eqnarray*}
}
It should also be marked that the limit of optimal stopping times as $K\to\infty$ coincide with the overall
optimal stopping time.

\section{\label{finalrem}Final remarks}
The presented results generalize known solutions of the optimal
stopping problems for the risk process (see \cite{fersie97:risk}, 
\cite{muc02:interest}) 
to more realistic models of
risk reserve processes. For  $\alpha=\beta=0$ the solution presented by \cite{fersie97:risk} 
are obtained. The model considered by \cite{muc02:interest} 
is not direct consequence of model (\ref{fer1}).
It implies that the solution of the optimal stopping problem for the risk reserve process investigated in
\cite{muc02:interest} 
is not simple conclusion from the formulae of Corollary \ref{fertw2}  and Section \ref{infcase}.

There are also similar optimal stopping problems considered by \cite{yas84:var} and  \cite{sch98:risk}.
The solution of the problem (\ref{fer6a}) is not direct consequence of the results from
neither \cite{yas84:var} nor \cite{sch98:risk}.


\end{document}